\documentclass{amsart}
\usepackage{latexsym,amsmath,amssymb,graphics,graphicx,amscd,amsfonts}
\usepackage[arrow, matrix, curve]{xy}
\usepackage[DIV11]{typearea}
\usepackage{amsthm}
\usepackage{eucal}
\usepackage{dsfont}

\newcommand{\C}{{\mathcal{C}}}
\newcommand{\D}{{\mathcal{D}}}
\newcommand{\V}{{\mathcal{V}}}
\newcommand{\K}{{\mathcal{K}}}
\newcommand{\T}{{\mathcal{T}}}
\newcommand{\I}{{\mathcal{I}}}
\newcommand{\R}{\mathbb{R}}

\newcommand{\mx}{\mathfrak{X}}
\newcommand{\lie}[1]{\mathfrak{#1}}
\newcommand{\dr}{\mathbf{d}}

\newcommand{\tg}{{\mathsf{t}}}
\newcommand{\s}{{\mathsf{s}}}
\newcommand{\m}{{\mathsf{m}}}
\newcommand{\rr}{{\rightrightarrows}}

\newcommand{\pair}[1]{\langle {#1}\rangle}

\DeclareMathOperator{\dom}{Dom}
\DeclareMathOperator{\erz}{span}
\allowdisplaybreaks[1]

\newtheorem{proposition}{Proposition}
\newtheorem{theorem}{Theorem}

\newtheorem{remark}{Remark}

\theoremstyle{definition}
\newtheorem{example}{Example}

\begin{document}
\title{Invariant generators for generalized distributions}
\author{M. Jotz}
\address{\textbf{Madeleine Jotz}\\ Section de Math{\'e}matiques\\ 
Ecole Polytechnique
  F{\'e}d{\'e}rale de Lausanne\\ 
1015 Lausanne\\ Switzerland
}
\email{madeleine.jotz@epfl.ch}
\author{T.S. Ratiu}
\address{ \textbf{Tudor S. Ratiu}\\ Section de Math{\'e}matiques\\
 et Centre Bernouilli\\ 
Ecole Polytechnique
  F{\'e}d{\'e}rale de Lausanne\\ 
1015 Lausanne\\ Switzerland}
\email{tudor.ratiu@epfl.ch}
\thanks{Partially supported by Swiss NSF grant 200021-121512.}
\date{}

\begin{abstract}
The existence of invariant
generators for distributions satisfying a
compatibility condition with the symmetry algebra is proved. 
\end{abstract}

\maketitle
\tableofcontents
\section{Introduction}

If a smooth manifold $M$ is acted upon in a proper way by a Lie group 
$G$, the space of orbits $M/G$ has the structure of a stratified space. 
If the action is free or with conjugated isotropy subgroups, 
the quotient $M/G$ is known 
to be a smooth manifold and the quotient map $\pi:M\to M/G$ a smooth surjective 
submersion. In the first case, a free and proper action is induced on the
tangent space $TM$ and on the cotangent space $T^*M$, 
but in the case of conjugated isotropies, the isotropy subgroups of the induced action
on the tangent space are not necessarily conjugated (see \cite{Rodriguez06} for 
a complete characterization of the isotropy lattice of the lifted action).

More generally, 
if a subdistribution of the Pontryagin bundle $\mathsf P_M:=TM\oplus T^*M$ is 
invariant under an action by an involutive subbundle $\mathcal I$
 of $TM$ how can we decide if  this 
distribution has invariant sections?  
The existence of invariant generators would imply that the generalized 
distribution, assumed that its cotangent part annihilates the vertical space of 
the action, pushes forward to a smooth generalized distribution on the space of leaves
$M/\mathcal I$. 
In this note,  we present a  theorem giving  sufficient conditions for a 
locally finitely generated
generalized distribution to be spanned by sections pushing forward to the quotient.  
This property has its origins in control theory and the first results in this 
direction were obtained in \cite{NivdS90} and \cite{ChTa89}.
Results about invariant generators for a certain class of invariant subbundles of exact Courant algebroids
are shown to hold in \cite{JoRaZa11}, but the techniques used there can not be applied 
in a straightforward manner to the present  situation.

\medskip

The paper is organized as follows. 
Background on smooth generalized distributions is rewieved in Section
\ref{two}.
Special emphasis is given to the notion of
pointwise and smooth annihilators since these are essential tools for the rest
of the paper. 
The main theorem of the paper,
guaranteeing, under certain assumptions, the 
existence of generators of generalized distributions that push forward to 
the quotient, is proved in Section
\ref{four}. 
\medskip

\paragraph{\textbf{Conventions and notations}}
If $M $ is a smooth manifold,  $C^\infty(M)$ denotes the sheaf of \emph{local}
functions on $M$, that is, 
an element $f \in C^\infty(M)$ is, by definition,  a smooth function $f:U\to
\R$, where the domain of definition $U$  of $f $ is an 
open subset of $M$. 
Similarly, if $E$ is a vector bundle over $M$, or a generalized
distribution on $M$, $\Gamma(E)$ denotes the set of smooth \emph{local}
sections of $E$. In particular, the sets of smooth  local
vector fields and one-forms on $M$ are denoted by  $\mx(M)$ and
$\Omega^1(M)$, respectively. The open domain of
definition of the local section $\sigma$ of $E$ is denoted by $\dom(\sigma)$.

\section{Generalities on distributions and smooth annihilators}\label{two}

\subsection{Smooth and pointwise annihilators}\label{dis}

The \emph{Pontryagin bundle} $\mathsf{P}:=TM \oplus T^* M$ of  a smooth manifold $M$ is naturally
endowed with a non-degenerate symmetric fiberwise bilinear form of signature
$(\dim M, \dim M)$ given by
\begin{equation}\label{pairing}
\left\langle (u_m, \alpha_m), ( v_m, \beta_m ) \right\rangle 
: = \left\langle\beta_m , u_ m \right\rangle + \left\langle\alpha_m, v _m \right\rangle
\end{equation}
for all $u _m, v _m \in T _mM$ and $\alpha_m, \beta_m \in T^\ast_mM$.

A  \emph{generalized distribution}  is a subset
$\Delta \subseteq \mathsf{P}=TM\oplus T^*M$  such that for each $m\in M$, the set 
$ \Delta (m) : = \Delta \cap \mathsf{P}(m)$ is a
vector subspace of $\mathsf{P}(m)= T_m M \oplus T_m ^\ast M$. The number $\dim \Delta
(m) $ 
is called the \emph{rank}  of $ \Delta $ at $ m  \in  M $. A point $m\in M$ is a
\emph{regular} point of the distribution $\Delta$ if there exists a
neighborhood $U$ of $m$ such that the rank of $\Delta$ is constant on
$U$. Otherwise, $m$ is a \emph{singular} point of the distribution. 

A local \emph{differentiable section} of $ \Delta $ is a smooth section 
$\sigma\in \Gamma(\mathsf{P})=\mx(M)\times\Omega^1(M)$
defined on some open subset $ U \subset  M $ such that  
$\sigma(u) \in  \Delta (u)$ for each $ u \in  U $;  
$ \Gamma(\Delta )$ denotes the space
of local differentiable sections of $ \Delta $. A generalized
distribution is said to be \emph{differentiable} or \emph{smooth} if for every
point $ m \in  M $ and every $ v \in  \Delta (m) $, there is a
differentiable section $ \sigma \in  \Gamma ( \Delta ) $ defined on an open
neighborhood $ U $ of $ m $ such that $ \sigma(m) = v $.

If $\Delta \subset \mathsf{P} $ is a smooth distribution, its
\emph{smooth orthogonal} distribution (or simply its \emph{smooth orthogonal}) is the smooth
generalized distribution $\Delta^\perp \subseteq \mathsf{P}$  defined by 
\begin{align*}
\Delta^ \perp (m): =
 \left\{(X(m),\alpha(m)) \left|  
\begin{array}{c}(X,\alpha) \in
\mx(M)\times\Omega^1(M) \text{ with } m\in\\ 
\dom(X)\cap\dom(\alpha) \text{ such that  for all }\\
(Y,\beta) \in\mx(M)\times\Omega^1(M)\\ 
\text{  with } m\in\dom(Y)\cap\dom(\beta),\\
\text{ we have }  \left\langle (X,\alpha),(Y,\beta)
\right\rangle = 0\\
 \text{ on } \dom(X)\cap\dom(Y)\cap\dom(\alpha)\cap\dom(\beta)
\end{array}\right.\right\}.
\end{align*}
In general, the inclusion $\Delta \subset \Delta^{ \perp\perp }$ is strict.
The smooth orthogonal  of a smooth generalized
 distribution
 is smooth by construction. If the distribution $\Delta$ is
a vector subbundle of $\mathsf{P}$, then its smooth orthogonal distribution is
also a vector subbundle of $\mathsf{P}$. Note also that the
smooth orthogonal of a smooth generalized distribution $\Delta$ is, in
general, different from the \emph{pointwise}
orthogonal distribution of $\Delta$, defined by
\begin{align*}
\Delta^{\perp_{\rm p}} (m): = \{\sigma_m\in \mathsf{P}(m)\mid
\pair{\sigma_m,\tau_m}=0 \text{ for all } \tau_m\in\Delta(m)\},
\end{align*}
where the subscript $p$ stands for ``pointwise''. The pointwise orthogonal
 of a smooth generalized distribution $\Delta$ is not smooth, in
general. The following statement is an immediate consequence of the definitions.

\begin{proposition}\label{point_smooth_annihi}
Let $\Delta$ be a smooth generalized distribution. Then we have
\[\Delta^\perp\subseteq \Delta^{\perp_{\rm p}},
\quad \Delta=\Delta^{\perp_{\rm p}\perp_{\rm p}},
\text{ and } \quad \Delta\subseteq\Delta^{\perp\perp}.\]
If $\Delta$ is itself a vector bundle over $M$, the smooth orthogonal
distribution $\Delta^\perp$ of
$\Delta$ is also a vector subbundle of $\mathsf{P}$, and  we have
$\Delta^\perp=\Delta^{\perp_{\rm p}}$.
\end{proposition}

This implies the following property of the smooth annihilator of a sum
of vector subbundles of $\mathsf{P}$; its proof is easy and can be found in \cite{JoRaSn11}.

\begin{proposition}\label{annihilators}
Let $\Delta_1$ and $\Delta_2$ be vector subbundles of the vector bundle
$(\mathsf{P},\langle\,,\rangle)$.
Since $\Delta_1$ and $\Delta_2$ have constant rank on $M$, their smooth
orthogonals $\Delta_1^\perp$ 
and $\Delta_2^\perp$ are also vector subbundles of $\mathsf{P}$ and equal to the
pointwise orthogonals of $\Delta_1$ and $\Delta_2$. The
intersection
$\Delta_1^\perp\cap\Delta_2^\perp$ is smooth if and only if  
it has locally constant rank.
\end{proposition}

\medskip

A tangent (respectively cotangent) distribution $\T\subseteq TM$ (respectively
$\mathcal{C}\subseteq T^*M$) can be identified with
the smooth generalized distribution $\T\oplus\{0\}$ 
(respectively $\{0\}\oplus\mathcal{C}$). 
The smooth orthogonal distribution of $\T\oplus\{0\}$ in $TM\oplus T^*M$ is easily computed to be 
$(\T\oplus\{0\})^\perp=TM\oplus \T^\circ$,
where
\[\T^\circ(m)=\left\{\alpha(m)\left|
\begin{array}{c}
 \alpha\in\Omega^1(M), m\in\dom(\alpha) \text{
  and } 
\alpha(X)=0\\ \text{ on }\dom(\alpha)\cap\dom(X)\text{ for all }
X\in\Gamma(\T)
\end{array}
\right.\right\}
\]
for all $m\in M$. This smooth cotangent distribution is called the smooth
\emph{annihilator} of $\T$. Analogously, we define the smooth annihilator
$\C^\circ$ of a
cotangent distribution $\C$. Then $\C^\circ$ is a smooth tangent distribution
and we have $(\{0\}\oplus\C)^\perp=\C^\circ\oplus T^*M$.

\subsection{Examples}
\begin{example}
If a Lie group $G$ with Lie algebra $\lie g$ acts on the manifold $M$, the
tangent distribution $\mathcal{V}$ whose value at each point $m \in M $ is given by
\[\V(m)=\{\xi_M(m)\mid \xi\in \mathfrak{g}\}\]
is called the \emph{vertical distribution}. 
Let $(\xi^1,\ldots,\xi^d)$ be a basis of the Lie algebra $\lie g$.
The smooth distribution $\V$
 is  tangent and globally finitely generated by the fundamental vector fields $\xi^1_M,\ldots,\xi^d_M$, but not of constant
rank, in general, unless the action is free or with conjugated isotropy subgroups.

If the $G $-action has conjugated isotropy groups, the
vertical distribution is a vector subbundle of $TM$. The smooth annihilator
$\V^\circ$ of $\V$ is given by
\begin{equation}\label{defK}
\V^\circ(m)
=\{\alpha(m)\mid \alpha\in\Omega^1(M), m\in\dom(\alpha),\text{ such that }
\alpha(\xi_M)=0 \text{ for all }\xi \in \lie g\}.
\end{equation}
In the following, we will also need the
 smooth generalized distribution $\K:=\V\oplus \{0\}$ and its smooth 
orthogonal $\K^\perp=TM\oplus\V^\circ$.
\end{example}

\begin{example}
Consider the tangent distribution $\mathcal D\subseteq T\R^2$ 
defined at every $(x,y) \in \mathbb{R}^2$ by
$$\mathcal D(x,y)=\left\{\begin{array}{lc}
\erz(\partial_y,\partial_x)& x\geq 0\\
\erz(\partial_y)& x<0
\end{array}\right..$$
Assume that $X$ is a vector field 
on $\R^2$ such that $X(0,0)=\partial_x$ and $X(x,y)\in\mathcal D(x,y)$
for all $(x,y)\in\R^2$. 
Write  $X=a\partial_x+b\partial_y$
with $a,b\in C^\infty(\R^2)$. Then $a(0,0)=1$ and $b(0,0)=0$.
Since $a$ is smooth, there exists a neighborhood $U$ of $(0,0)$ in $\R^2$ such that
$a$ doesn't vanish on $U$. This implies that $X$ does not
take values in $\mathcal D$, a contradiction. Hence, $\mathcal D$ is not a smooth distribution.
\end{example}

\begin{example}\label{smooth_dis}
Consider the tangent distribution $\mathcal D\subseteq T\R^2$ 
defined at every $(x,y) \in \mathbb{R}^2$ by
$$\mathcal D(x,y)=\left\{\begin{array}{lc}
\erz(\partial_y,\partial_x)& x> 0\\
\erz(\partial_y)& x\leq 0
\end{array}\right..$$
This distribution is smooth.
\end{example}

\section{Invariant generators for distributions}\label{four}
\subsection{The theorem}
We present here a theorem which can help to decide if a locally
finitely generated distribution is spanned by its descending sections. 
The proof is inspired by \cite{ChTa89}.
  
The space $\Gamma(TM \oplus T ^\ast M) $ of local sections of the Pontryagin
bundle is endowed with a skew-symmetric bracket given by
\begin{align}\label{wrong_bracket}
[(X, \alpha), (Y, \beta) ] : 
&= \left( [X, Y],  \boldsymbol{\pounds}_{X} \beta 
- \boldsymbol{\pounds}_{Y} \alpha + \frac{1}{2} \mathbf{d}\left(\alpha(Y) 
- \beta(X) \right) \right) \nonumber \\
&= \left([X, Y],  \boldsymbol{\pounds}_{X} \beta - \mathbf{i}_Y
  \mathbf{d}\alpha - \frac{1}{2} \mathbf{d} \left\langle (X, \alpha), (Y, \beta) \right\rangle
\right)
\end{align}
(see \cite{Courant90a}). This bracket is $\mathds{R}$-bilinear (in the sense that 
$[a_1(X_1, \alpha_1)+a_2(X_2, \alpha_2), (Y, \beta) ]
=a_1[(X_1, \alpha_1), (Y, \beta) ]+a_2[(X_2, \alpha_2), (Y, \beta)
]$ for all $a_1,a_2\in\R$ and
$(X_1, \alpha_1),(X_2, \alpha_2), (Y, \beta)\in\Gamma(TM\oplus T^*M)$  on the common domain of
definition of the three sections) and does not satisfy the Jacobi identity.

Let $\mathcal{I} \subseteq  TM$ be a smooth tangent distribution. If $(Y,0)$
is a local section of  $\I\oplus\{0\}\subseteq TM\oplus
T^*M$ and $(X,\alpha)$ a local section of  $TM\oplus \I^\circ$, using  
$\alpha(Y)=0$ we get 
\begin{equation}
\label{I_bracket} 
[(Y, 0), (X, \alpha) ] = \left( [Y, X],  \boldsymbol{\pounds}_{Y} \alpha  
+ \frac{1}{2} \mathbf{d}\left( - \alpha(Y) \right) \right)
=\left( \boldsymbol{\pounds}_{Y} X,  \boldsymbol{\pounds}_{Y} \alpha \right).
\end{equation} 
Note also that
\begin{equation}
\label{I_bracket_function} 
[(Y, 0), f(X, \alpha) ]  = \left( \boldsymbol{\pounds}_{Y} f \right) 
\left( X, \alpha \right) + f[(Y, 0), (X, \alpha) ] 
\end{equation} 
for all $f\in C^\infty(M)$.

\begin{theorem}\label{prop-ch-ta}
Let $\mathcal{I}  \subseteq  TM$ be  an involutive vector subbundle of  
 $TM$ and $\mathcal{D}\subseteq TM\oplus \I^\circ$ a 
generalized distribution on $M$. Let $\Theta:=\I\oplus\{0\}\subseteq TM\oplus T^*M$. Assume  that for
  each   $m \in M$  there exist an open set $ U \subseteq M$ with $ m \in  U $
and smooth sections $d_1,\ldots,d_r\in\Gamma(\mathcal D)$
such that 
\[\mathcal D(m')=\erz_\R\{d_1(m'), \ldots, d_r(m')\}\]
for all $m'\in U$ and 
\begin{equation}\label{ass1}
[d_i,\Gamma(\Theta)]\subseteq\erz_{C^\infty(U)}\{d_1,\ldots,d_r\}+\Gamma(\Theta)
\end{equation}
for $i=1,\ldots,r$.
Let $(X,\alpha)$ be a local section of $TM\oplus \I^\circ$ satisfying 
\begin{equation}\label{ass2}
[(X, \alpha),\Gamma(\Theta)]\subseteq\erz_{C^\infty(U)}\{d_1,\ldots,d_r\}+\Gamma(\Theta).
\end{equation}
Then there exists 
smooth sections $d$, $d_1', \ldots,d_r'$ in  $
\erz_{C^\infty(U)}\{d_1,\ldots,d_r\}$ satisfying
\begin{enumerate}
\item[{\rm (i)}]
  $\mathcal{D}(m')=\operatorname{span}\{d_1'(m'),\ldots,d_r(m')\}$ for
all $m'\in U$,
\item[{\rm (ii)}] $[d_i',\Gamma(\Theta)]\subseteq\Gamma(\Theta)$ on $U$ for
all $i=1,\dots,r$, and
\item[{\rm (iii)}] $[(X,\alpha)+d,\Gamma(\Theta)]\subseteq\Gamma(\Theta)$ on $U$.
\end{enumerate}
\end{theorem}

\begin{remark}{\rm 
Note that if $\mathcal D$ has constant rank on $U$, then
$\Gamma(\mathcal D)=\erz_{C^\infty(U)}\{d_1,\ldots,d_r\}$. For
a generalized distribution this is not necessarily true.}
\end{remark}

\subsection{The proof}
This is a long proof and so it will be broken up in four steps. Since the
statement is local, we work in a foliated chart of $ \mathcal{I}$. We choose a
set of local sections as in the hypothesis of the theorem
and write its elements as a
sum of  a component tangent to the leaves of $\mathcal{I}$ and the rest. The
main work is the analysis of this second component. The $r$ local sections
spanning pointwise $\D$
in the theorem are constructed from an initially chosen  set of local spanning
sections of $\mathcal{D}$ using in an essential way the information gathered
about their second component. In the first step we construct a family of
linear systems for the derivatives (along the leaves of $\I$) of these second components.  
Using specific properties of this system, in the second step, we linearly
transform the second components in order to get pairs formed by a vector field
and a one-form that are independent of the coordinates of the leaves of $
\mathcal{I} $.  In the third step, we extend this linear transformation to 
the set of local spanning sections of $ \mathcal{D} $ and, using the property 
found in the previous step, $r$ sections of $\mathcal{D} $ are constructed
that satisfy the first two properties in the statement. In the fourth step an 
additional section of $ \mathcal{D} $  is constructed that satisfies the third 
property in the statement.

 \medskip

\paragraph{\textit{Step 1: Construction of the linear system.}}
Let $n: = \dim M $ and $k: = \dim \I(x)$, for $ x \in  M $. Since the vector
subbundle $ \I $ is involutive, it is integrable by the Frobenius Theorem
and thus any $ m \in  M $ lies in a foliated chart domain $ U _1 $ described by
coordinates $(x^1,\dots,x^n)$ such that the first $ k $ among them define the 
local integral submanifold containing $ m $. Thus, for any $ m' \in  U _1 $
the basis vector fields $\partial_{x^1},\dots,\partial_{x^k} $ evaluated at $ m' $ 
span $\mathcal{I}(m') $.

Because $\mathcal{D}$ is locally finitely generated, we can find on a sufficiently small
neighborhood $U$ of $m$ in $U_1$, smooth  sections
$d_1=(X_1,\alpha^1),\dots,d_r=(X_r,\alpha^r)$ spanning $\Gamma_U(\mathcal{D})$ as a $C ^{\infty}(U)$-module.
Write, for $i=1,\dots,r$,
\[
d_i=(X_i,\alpha^i)=\sum_{j=1}^n(X_i^j\partial_{x^j},\alpha^i_j\dr{x^j}) ,
\]
with $ X_i^j$ and $\alpha_j^i$ smooth local functions defined on $U$ for
$j=1,\dots,n$. Note that  $\alpha_1^i=\ldots=\alpha_k^i=0$ for $i=1,\dots, r$, 
since $\alpha^i\in\Gamma(\I^\circ)$.
By hypothesis \eqref{ass1} and with $\partial_{x^l}\in\Gamma(\I)$ for
$l=1,\dots, k$, we get for all $i=1,\dots,r$ and $l=1,\dots,k$:
\begin{align*}
\partial_{x^l}(X_i,\alpha^i)&:=\left[\left(\partial_{x^l},0 \right),\left(X_i,\alpha^i\right)\right]
=\left(\left[\partial_{x^l},X_i\right],\boldsymbol{\pounds}_{\partial_{x^l}}\alpha^i\right)\\
&=\sum_{j=1}^n\left(\partial_{x^l}(X_i^j)\partial_{x^j},\partial_{x^l}(\alpha^i_j)\mathbf{d}x^j\right)
\in \erz_{C^\infty(U)}\{d_1,\ldots,d_r\}+\Gamma(\Theta).
\end{align*}
Hence we can write
\[\partial_{x^l}(X_i,\alpha^i)
=\sum_{j=1}^n\left(\partial_{x^l}(X_i^j)\partial_{x^j},\partial_{x^l}(\alpha_j^i)\mathbf{d}x^j\right)
=\sum_{j=1}^kA_{li}^j(\partial_{x^j},0)+\sum_{s=1}^rB_{ls}^i\left(X_s,\alpha^s\right)
\]
with  $ A_{li}^j ,
B_{ls}^{i}\in C^\infty(U)$ for $i,s=1,\dots,r$ and $l, j=1,\dots,k$.
Setting 
\[\left(\widetilde{X}_i,\widetilde{\alpha}^i\right)
:=\sum_{j=k+1}^n\left(X_i^j\partial_{x^j},\alpha_j^i\dr x^j \right)
\] 
for $i=1,\dots,r$, we get
\begin{align}
\label{tail_identity}
\partial_{x^l}\left(\widetilde{X}_i,\widetilde{\alpha}^i\right)
:&=\left[ \Big(\partial_{x^l},0\Big),\left(\widetilde{X}_i,\widetilde{\alpha}^i\right)\right]
=\sum_{j=k+1}^n\left(\partial_{x^l}(X_i^j)\partial_{x^j},\partial_{x^l}(\alpha_j^i)\dr
x^j\right)  \\
&=\sum_{s=1}^rB_{ls}^i\left(\widetilde{X}_s,\widetilde{\alpha}^s\right). \nonumber 
\end{align}
We verify the last equality. Since $ \alpha^s_j = 0 $ for $ j= 1, \ldots, k$ 
and $ s = 1, \ldots, r $, we have for any $ i = 1, \ldots, r$, $ l = 1, \ldots, k$,
\begin{align*}
\partial _{x^l} &(X _i , \alpha ^i ) 
 = \sum_{m=1}^k\left(\partial_{x^l}(X_i^m)\partial_{x^m},0\right) + 
\sum_{j=k+1}^n\left(\partial_{x^l}(X_i^j)\partial_{x^j},\partial_{x^l}(\alpha_j^i)\mathbf{d}x^j\right) \\
& = \sum_{m=1}^kA_{li}^m(\partial_{x^m},0)+\sum_{s=1}^rB_{ls}^i\left(X_s,\alpha^s\right) \\
& =  \sum_{m=1}^kA_{li}^m(\partial_{x^m},0) +
\sum_{s=1}^rB_{ls}^i \sum_{j=1}^n(X_s^j\partial_{x^j},\alpha^s_j\dr{x^j}) \\
& =  \sum_{m=1}^kA_{li}^m(\partial_{x^m},0) +
\sum_{s=1}^rB_{ls}^i \sum_{m=1}^k(X_s^m\partial_{x^m},0)
+ \sum_{s=1}^rB_{ls}^i \sum_{j=k+1}^n(X_s^j\partial_{x^j},\alpha^s_j\dr{x^j}) \\
& = \sum_{m=1}^k\left(A_{li}^m + \sum_{s=1}^rB_{ls}^i X_s^m \right) (\partial_{x^m},0)
+ \sum_{j=k+1}^n \left(\left( \sum_{s=1}^rB_{ls}^i X_s ^j \right) \partial
  _{x^j} , \left( \sum_{s=1}^rB_{ls}^i \alpha^s_j \right) \mathbf{d} x^j  \right) 
\end{align*} 
which is equivalent to
\begin{align} 
\label{first_tail_identity}
\partial _{x ^l }(X^m_i) &= A_{li}^m + \sum_{s=1}^rB_{ls}^i X_s^m, \quad
\text{for all } \quad  l,m = 1, \ldots, k,\: \: i=1, \ldots, r \\
\label{second_tail_identity}
\partial _{x ^l }(X^j_i) &= \sum_{s=1}^rB_{ls}^i X_s^j \quad \text{for all }
\quad l = 1, \ldots, k, \:\:  i=1, \ldots, r, \: \: j = k+1, \ldots n\\
\label{third_tail_identity}
\partial _{x^l}( \alpha ^i _j) & = \sum_{s=1}^rB_{ls}^i \alpha^s_j  
\quad \text{for all } \quad l = 1, \ldots, k, \:\:   i=1, \ldots, r, \:\:  j = k+1, \ldots n.
\end{align} 
Using \eqref{second_tail_identity} and \eqref{third_tail_identity} we get
\begin{align*}
\partial_{x^l}\left(\widetilde{X}_i,\widetilde{\alpha}^i\right)
&=\sum_{j=k+1}^n\left(\partial_{x^l}(X_i^j)\partial_{x^j},\partial_{x^l}(\alpha_j^i)\mathbf{d} x^j\right) \\
& = \sum_{j=k+1}^n\left(\left(\sum_{s=1}^rB_{ls}^i X_s^j \right) \partial_{x^j}, 
\left( \sum_{s=1}^rB_{ls}^i \alpha^s_j \right) \mathbf{d} x^j\right) \\
& = \sum_{j=k+1}^n \sum_{s=1}^rB_{ls}^i \left( X_s^j \partial_{x^j}, 
\alpha^s_j \mathbf{d} x^j\right) 
 = \sum_{s=1}^rB_{ls}^i \sum_{j=k+1}^n \left( X_s^j \partial_{x^j}, 
\alpha^s_j \mathbf{d} x^j\right) \\
&=\sum_{s=1}^rB_{ls}^i\left(\widetilde{X}_s,\widetilde{\alpha}^s\right)
\end{align*} 
which proves \eqref{tail_identity}. Note that \eqref{first_tail_identity} gives 
a formula for $ \partial _{x^l}(X^m_i)$ for all $l,m = 1, \ldots, k$ and 
$ i=1, \ldots, r$.
\medskip

\paragraph{\textit{Step2: Construction of an $r\times r$-matrix $B$ such that $
\left((\widetilde{X}_{1},\widetilde{\alpha}^1),\dots,
(\widetilde{X}_{r},\widetilde{\alpha}^r)\right) B$ does not depend on $x^1,\dots,x^k$.}}
We rewrite the system \eqref{tail_identity} in the form
\begin{equation}
\label{tail_system}
\left( \partial_{x^l}\left(\widetilde{X}_{1},\widetilde{\alpha}^1\right),
  \ldots, \partial_{x^l}\left(\widetilde{X}_r,\widetilde{\alpha}^r\right)\right) =
\left(\left(\widetilde{X}_{1},\widetilde{\alpha}^1\right) \ldots,
\left(\widetilde{X}_r,\widetilde{\alpha}^r\right)\right) B_l, 
\end{equation} 
where $ B_l := [B^i_{ls}]$ is the $r \times r$ matrix whose entry $
B_{ls}^{i} \in  C ^{\infty}(U)$ is the intersection of the $i$-th column and
the $s$-th row ($ i, j =1, \ldots, r$). In view of
\eqref{second_tail_identity}, \eqref{third_tail_identity}, this system can be explicitly written as
\begin{equation*} 
\partial_{x^l} \left[ 
\begin{array}{ccc}
X_1^{k+1}& \ldots &X_r^{k+1}\\
\vdots& \vdots & \vdots\\
X_1^{n}& \ldots &X_r^{n}\\
\alpha^1_{k+1}& \ldots & \alpha^r_{k+1}\\
\vdots& \vdots & \vdots\\
\alpha^1_{n}& \ldots & \alpha^r_{n}\\
\end{array} \right] = 
\left[ 
\begin{array}{ccc}
X_1^{k+1}& \ldots &X_r^{k+1}\\
\vdots& \vdots & \vdots\\
X_1^{n}& \ldots &X_r^{n}\\
\alpha^1_{k+1}& \ldots & \alpha^r_{k+1}\\
\vdots& \vdots & \vdots\\
\alpha^1_{n}& \ldots & \alpha^r_{n}\\
\end{array} \right]
\left[ 
\begin{array}{cccc}
B^1_{l1} & B^2_{l1} & \ldots & B^r_{l1}\\
B^1_{l2} & B^2_{l2} & \ldots & B^r_{l2}\\
\vdots&\vdots&\vdots&\vdots\\ 
B^1_{lr} & B^2_{lr} & \ldots & B^r_{lr}
\end{array} 
\right] ,
\end{equation*} 
where 
\begin{equation}
\label{matrix_formula_x_alpha}
\left((\widetilde{X}_{1},\widetilde{\alpha}^1),\dots,(\widetilde{X}_{r},\widetilde{\alpha}^r)\right)=
\left[ 
\begin{array}{ccc}
X_1^{k+1}& \ldots &X_r^{k+1}\\
\vdots& \vdots & \vdots\\
X_1^{n}& \ldots &X_r^{n}\\
\alpha^1_{k+1}& \ldots & \alpha^r_{k+1}\\
\vdots& \vdots & \vdots\\
\alpha^1_{n}& \ldots & \alpha^r_{n}\\
\end{array} \right].
\end{equation} 
Equivalently, taking the transpose of this system, we get
\begin{align}
\label{transpose_system}
&\partial_{x^l} \left[ 
\begin{array}{cccccc}
X_1^{k+1}& \ldots &X_1^n & \alpha _{k+1}^1& \ldots &\alpha _n ^1\\
X_2^{k+1}& \ldots &X_2^n & \alpha _{k+1}^2& \ldots &\alpha _n ^2\\
\vdots & \vdots & \vdots & \vdots & \vdots & \vdots \\
X_r^{k+1}& \ldots &X_r^n & \alpha _{k+1}^r& \ldots &\alpha _n ^r
\end{array} 
\right] = B_l ^\top
\left[ 
\begin{array}{cccccc}
X_1^{k+1}& \ldots &X_1^n & \alpha _{k+1}^1& \ldots &\alpha _n ^1\\
X_2^{k+1}& \ldots &X_2^n & \alpha _{k+1}^2& \ldots &\alpha _n ^2\\
\vdots & \vdots & \vdots & \vdots & \vdots & \vdots \\
X_r^{k+1}& \ldots &X_r^n & \alpha _{k+1}^r& \ldots &\alpha _n ^r
\end{array}
\right]
\end{align} 

Now fix $j\in\{1,\dots,k\}$, think of $ x^j $ as a time variable and all the 
other $x^i$, $i \neq j $,  as parameters. Consider the following linear 
ordinary differential equation on $ \mathbb{R} ^r $, where the solutions $Y$ depend smoothly on the parameters $x^1, \ldots, x^{j-1}, x^{j+1}, \ldots , x^n$:
\begin{equation}\label{diffeq}
\partial _{x^j} Y = B_j^\top Y. 
\end{equation} 
Let $Y^j_1,\dots,Y^j_r$ be $r $ linearly independent solutions of
\eqref{diffeq}. The $ r \times  r $ matrix $ W _j $ whose columns are
$Y^j_1,\dots,Y^j_r$, that is, $W_j: =(Y^j_1,\dots,Y^j_{r})$, is invertible. 
However, by \eqref{transpose_system} with $l = j $, the $ 2(n-k)$ columns of 
the matrix in this system also satisfy \eqref{diffeq} and therefore these 
columns are linear combinations of $Y^j_1,\dots,Y^j_r$, that is, there exists a matrix 
$L_j$ having $ r $ rows and $2(n-k) $ columns such that 
\begin{equation*}
\left((\widetilde{X}_{1},\widetilde{\alpha}^1),\dots,(\widetilde{X}_{r},\widetilde{\alpha}^r)\right)^\top
=W_jL_j,\quad j=1,\dots,k.
\end{equation*}
The entries of $ L_j $ are smooth functions of $(x^1, \ldots, x^n)$ and are
independent of the variable $x^j$ (the ``time'' in the differential equation \eqref{diffeq}).

Since this argument holds for any $ j=1, \ldots, k $ this implies 
\begin{equation}\label{idwl}
W_1L_1=W_2L_2=\dots=W_kL_k.
\end{equation}
Because $W_2$ is nonsingular, we have 
\[L_2(x^1,\ldots,x^n)=\left(W_2^{-1}W_1L_1\right)(x^1,\ldots,x^n).\]
Since $L_2$ is independent of $x^2$ and $L_1$ is independent of $x^1$, we get
\begin{align}
\label{L_2}
L_2(x^1,x^2, \ldots,x^n)&=L_2(x^1,0,x^3,\ldots,x^n) \\
&=\left(W_2^{-1}W_1\right)(x^1,0,x^3,\ldots,x^n)L_1(x^1,0,x^3\ldots,x^n)
\nonumber \\
&=\left(W_2^{-1}W_1\right)(x^1,0,x^3,\ldots,x^n)L_1(0,0,x^3\ldots,x^n). \nonumber 
\end{align}
In the same manner, we have
$L_3(x^1,\ldots,x^n)=\left(W_3^{-1}W_2L_2\right)(x^1,\ldots,x^n)$,
and hence
\begin{align*}
L_3(x^1,\ldots,x^n)&=L_3(x^1,x^2,0,x^4,\ldots,x^n)\\
&=\left(W_3^{-1}W_2\right)(x^1,x^2,0,x^4,\ldots,x^n) L_2(x^1,x^2,0,x^4,\ldots,x^n)\\
&\overset{\eqref{L_2}}=\left(W_3^{-1}W_2\right)(x^1,x^2,0,x^4,\ldots,x^n)
\left(W_2^{-1}W_1\right)(x^1,0,0,x^4,\ldots,x^n)\\
&\quad\quad \times  L_1(0,0,0,x^4\ldots,x^n).
\end{align*}
By induction, we get
\begin{align}
\label{L_k}
L_k(x^1,&\ldots,x^n)=L_k(x^1,\ldots,x^{k-1},0,x^{k+1},\ldots,x^n)\\
&=\left(W_k^{-1}W_{k-1}\right)(x^1,\ldots,x^{k-1},0,x^{k+1},\ldots,x^n) \nonumber \\
& \qquad \times  L_{k-1}(x^1,\ldots,x^{k-1},0,x^{k+1},\ldots,x^n) \nonumber \\
&=\left(W_k^{-1}W_{k-1}\right)(x^1,\ldots,x^{k-1},0,x^{k+1},\ldots,x^n)
\nonumber \\
& \qquad \times   \left(W_{k-1}^{-1}W_{k-2}\right)(x^1,\ldots, x^{k-2},0, 0,
x^{k+1},\ldots,x^n) \times  \cdots \nonumber\\
& \qquad \times  \left(W_3^{-1}W_2\right)(x^1, x^{2},0,\ldots, 0, x^{k+1},\ldots,x^n)\nonumber\\
&\qquad \times  \left(W_2^{-1}W_1\right)(x^1,0,\ldots, 0,x^{k+1},\ldots,x^n)
L_1(0,\ldots,0,x^{k+1}\ldots,x^n). \nonumber 
\end{align}
Define the smooth $ r \times  r $ nonsingular matrix
\begin{align}
\label{H}
H(x^1,\ldots,x^n)&
:=W_k(x^1,\ldots,x^n)\left(W_k^{-1}W_{k-1}\right)(x^1,\ldots,x^{k-1},0,x^{k+1},\ldots,x^n)\\
&\qquad \times  \cdots \times  \left(W_3^{-1}W_2\right)(x^1,x^{2},0,\ldots, 0, x^{k+1},\ldots,x^n)\nonumber\\
&\qquad \times  \left(W_2^{-1}W_1\right)(x^1,0,\ldots, 0,x^{k+1},\ldots,x^n)
\nonumber 
\end{align}
 and the smooth $r\times (2n-2k)$ matrix
\begin{align*}
L(x^1,\ldots,x^n):= L_1(0,\ldots,0,x^{k+1},\ldots,x^n).
\end{align*}
Note that $L$ does not depend on $x^1,\ldots,x^k$.
Using equations \eqref{L_k} and \eqref{H}, we get
\begin{align*}
(HL)(x^{1},\ldots,x^n) &= H(x^1, \ldots, x^n) 
L_1(0, \ldots, 0, x^{k+1}, \ldots, x^n)\\
&=W_k(x^1,\ldots,x^n)L_k(x^1,\ldots,x^n)\\
&=\left((\widetilde{X}_{1},\widetilde{\alpha}^1),\dots,(\widetilde{X}_{r},\widetilde{\alpha}^r)\right)^\top.
\end{align*} 
Define the $r \times r$ nonsingular matrix $B$ depending smoothly on $(x^1, \ldots, x^n)$ by 
\begin{equation}\label{defb}
B(x^1, \ldots, x^n):=(H(x^1, \ldots, x^n)^\top)^{-1}=(H(x^1, \ldots, x^n)^{-1})^\top.
\end{equation}
Then
$\left((\widetilde{X}_{1},\widetilde{\alpha}^1),\dots,(\widetilde{X}_{r},\widetilde{\alpha}^r)\right)
B=L^\top$ does \textit{not} depend on $x^1,\dots,x^k$.  

\medskip

\paragraph{\textit{Step 3: Construction of the local sections $(Z_1, \gamma _1), \ldots, (Z_r, \gamma _r)$.}}
We want to better understand the columns of this matrix.
 In view of \eqref{matrix_formula_x_alpha}, we can write
\[
\left((\widetilde{X}_{1},\widetilde{\alpha}^1),\dots,(\widetilde{X}_{r},\widetilde{\alpha}^r)\right) B
 =\left[\begin{array}{c} \widetilde{C}\\
\widetilde{D}\end{array}\right],
 \]
where
\[ 
\widetilde{C}=\left[\widetilde{C}_{jl}\right]_{\substack{
j=k+1,\dots,n\\l=1,\dots,r}}\quad\text{ and }\quad
\widetilde{D}=\left[\widetilde{D}_{jl}\right]_{\substack{j=k+1,\dots,n\\l=1,\dots,r}}
\] 
are $(n-k)\times r$-matrices whose entries are smooth functions of only
$x^{k+1}, \ldots, x^n$ (that is, they do not depend on $x^1, \ldots, x^k$). Thus the 
$i$-th column of the matrix
$\left((\widetilde{X}_{1},\widetilde{\alpha}^1),\dots,(\widetilde{X}_{r},\widetilde{\alpha}^r)\right)
B$ equals
\begin{equation}\label{id}
\left(\left((\widetilde{X}_{1},\widetilde{\alpha}^1),\dots,(\widetilde{X}_{r},\widetilde{\alpha}^r)\right)
B\right)_i
:=\sum_{j=k+1}^n(\widetilde{C}_{ji}\partial_{x^j},\widetilde{D}_{ji}\dr x^j) 
\end{equation}
and so, by \eqref{I_bracket}, we get for any $l=1, \ldots, k$ and $i=1, \ldots, r$, 
\begin{equation}\label{diff}
\left[(\partial_{x^l},0),
\left(\left((\widetilde{X}_{1},\widetilde{\alpha}^1),\dots,(\widetilde{X}_{r},\widetilde{\alpha}^r)\right) 
B\right)_i\right]=0.
  \end{equation} 

Let $d_1'=(Z_1,\gamma_1),\dots,d_r'=(Z_r,\gamma_r)$ be the sections of 
$ \mathcal{D} \subseteq TM \oplus \mathcal{I} ^ \circ$  on $ U $ defined by
\[
d_i'=(Z_i,\gamma_i):=\Big(\big((X_{1},\alpha^1),\dots,(X_r,\alpha^r)\big) B\Big)_i
\quad \text{for} \quad i=1,\dots,r,
\]
where, as before, if  we write
\[
\big((X_{1},\alpha^1),\dots,(X_r,\alpha^r)\big) B
= \left[ 
\begin{array}{ccc}
X_1^{1}& \ldots &X_r^{1}\\
\vdots& \vdots & \vdots\\
X_1^{k}& \ldots &X_r^{k}\\
X_1^{k+1}& \ldots &X_r^{k+1}\\
\vdots& \vdots & \vdots\\
X_1^{n}& \ldots &X_r^{n}\\
\alpha^1_{1}& \ldots & \alpha^r_{1}\\
\vdots& \vdots & \vdots\\
\alpha^1_{k}& \ldots & \alpha^r_{k}\\
\alpha^1_{k+1}& \ldots & \alpha^r_{k+1}\\
\vdots& \vdots & \vdots\\
\alpha^1_{n}& \ldots & \alpha^r_{n}\\
\end{array} \right]
\left[ 
\begin{array}{cccc}
B_{11} & B_{12} & \ldots & B_{1r}\\
B_{21} & B_{22} & \ldots & B_{2r}\\
\vdots&\vdots&\vdots&\vdots\\ 
B_{r1} & B_{r2} & \ldots & B_{rr}
\end{array} 
\right] 
=\left[\begin{array}{c}C\\
D\end{array}\right],
 \]
with
\[ 
C=\big[C_{jl}\big]_{\substack{j=1,\dots,n\\l=1,\dots,r}}\quad\text{ and }\quad
D=\big[D_{jl}\big]_{\substack{j=1,\dots,n\\l=1,\dots,r}}\] 
$n\times r$-matrices whose entries depend smoothly on all coordinates $x^1, \ldots, x^n$, we have
\[
\Big(\big((X_{1},\alpha^1),\dots,(X_r,\alpha^r)\big) B\Big)_i
= \sum_{j=1}^n(C_{ji}\partial_{x^j},D_{ji}\dr x^j).
\]
Note that since $\alpha^i_j=0$ for $i=1,\dots,r$ and $j=1,\dots,k$, we get
$D_{ji}=0$ for all  $i=1,\dots,r$ and $j=1,\dots,k$.
Using \eqref{id}, we conclude
\begin{align*}
Z_i&=\xi_i+ \mathbb{P}_{TM} \left( \left( \left(
      \widetilde{X}_{1},\dots,\widetilde{X}_{r} \right)  
B\right)_i \right) 
=\xi_i+\sum_{j=k+1}^nC_{ji}\partial_{x^j} \,,\\
\gamma _i 
& = \mathbb{P}_{T^*M} \left( \left( \left(
      \widetilde{\alpha}_{1},\dots,\widetilde{\alpha}_{r}\right) 
B\right) _i \right)
= \sum_{j=k+1}^nD_{ji}\mathbf{d} x^j \,,
\end{align*} 
where $ \mathbb{P}_{TM}$ and $\mathbb{P}_{T^*M}$ are the projections 
on the vector field and one-form factors, respectively, and
\begin{align*}
\xi_i&=\sum_{j=1}^kC_{ji}\partial_{x^j}\in \Gamma(\I).
\end{align*}
Thus, we have for all $l=1,\dots,k$ and $i=1,\dots,r$,
\begin{align} 
\label{int_rel_first_part}
\left[\left(\partial_{x^l},0\right),\left(Z_i,\gamma_i\right)\right]
& =\left[\left(\partial_{x^l},0\right),\left(\xi_i,0\right)\right]
+\left[\left(\partial_{x^l},0\right),\left(\left((\widetilde{X}_{1},\widetilde{\alpha}_1),\dots,
(\widetilde{X}_{r},\widetilde{\alpha}_r)\right)\cdot B\right)_i\right]  \\
&\overset{\eqref{diff}}=([\partial_{x^l},\xi_i],0)+0\in\Gamma(\Theta) \nonumber 
\end{align}
since $ \mathcal{I} \subset  TM$ is an involutive vector subbundle, by hypothesis. 
Hence, if we write an arbitrary section $(\eta,0) \in \Gamma(\Theta)$ as
\[
(\eta,0)=\left(\sum_{j=1}^k\eta_j\partial_{x^j},0\right),
\]
where $\eta_1,\dots,\eta_k$ are smooth functions of $x^1, \ldots, x^n$, 
we get for $i=1,\dots,r$
\begin{align*}
\left[\left(\eta,0\right),\left(Z_i,\gamma_i\right) \right]
&=\left[\left(\sum_{j=1}^k\eta_j\partial_{x^j},0\right),(Z_i,\gamma_i)\right]\\
&=\sum_{j=1}^k\eta_j \left[\left(\partial_{x^j},0\right),\left(Z_i,\gamma_i\right) \right] 
-\sum_{j=1}^kZ_i[\eta_j](\partial_{x^j},0)\in\Gamma(\Theta).
\end{align*}
Indeed, both terms are elements of $\Gamma(\Theta)$:
the first summand by \eqref{int_rel_first_part} and the second summand because of its form.

Thus, since by construction, $(Z_1,\gamma_1),\dots,(Z_r,\gamma_r)$ span
$\mathcal{D}$ on $U$ (because $B$ is an invertible $ r \times  r $  matrix),
these smooth sections of $\mathcal{D} \subseteq TM\oplus\I^\circ$ satisfy the first two
statements of the proposition.

\medskip

\paragraph{\textit{Step4: Construction of the local section $d$.}}
We use the sections $d_1=(X_1, \alpha^1), \ldots, d_r=(X_r, \alpha^r)$ 
spanning $\mathcal{D}$ on the open subset $U \subseteq M$. 
Since $(X, \alpha) \in  \Gamma(TM \oplus \mathcal{I}^\circ)$, it can be written in the form
\[
(X,\alpha)=\sum_{j=1}^na^j(\partial_{x^j},0)+\sum_{j=k+1}^nb_j(0,\dr {x^j}),
\]
where $a^1,\dots,a^n, b_{k+1},\dots,b_n$ are $C^\infty$-functions of $x^1,
\ldots, x^n$. By hypothesis \eqref{ass2}, 
\[
\left[(\partial_{x^l},0), (X,\alpha)\right]
=\sum_{j=1}^n\partial_{x^l}(a^j)(\partial_{x^j},0)+\sum_{j=k+1}^n\partial_{x^l}(b_j)(0,\dr {x^j})\in
\erz_{C^\infty(U)}\{d_1,\ldots,d_r\}+\Gamma(\Theta)
\] 
for all $l=1,\dots,k$. Thus, for each $l=1, \ldots, k$, there exist functions
$\beta_l^1,\dots,\beta_l^r$ and $\sigma_l^1,\dots,\sigma_l^k$ depending
smoothly on $x^1, \ldots, x^n$ such that 
\begin{align*}
\sum_{j=1}^n\partial_{x^l}(a^j)(\partial_{x^j},0)+\sum_{j=k+1}^n\partial_{x^l}(b_j)(0,\dr {x^j})
&=\sum_{j=1}^k\sigma_l^j(\partial_{x^j},0)+\sum_{j=1}^r\beta_l^j(X_j,\alpha^j).
\end{align*}
Hence, if we define 
\[
(\widetilde{X},\tilde\alpha):=\sum_{j=k+1}^n(a^j\partial_{x^j},b_j\dr x^j),
\] 
then proceeding as in the proof of \eqref{tail_identity}, we get for each $l=1,\ldots,k$,
\begin{align} \label{bracket_tilde_l}
\left[(\partial_{x^l},0),(\widetilde{X},\tilde\alpha)\right]
&=\sum_{j=k+1}^n(\partial_{x^l}(a^j)\partial_{x^j},\partial_{x^l}(b_j)\dr{x^j})\\
&=\sum_{j=1}^r\beta_l^j(\widetilde{X}_j,\tilde\alpha^j)
=\left((\widetilde{X}_{1},\tilde\alpha^1),\dots,(\widetilde{X}_r,\tilde\alpha^r)\right)\beta_l \nonumber
\end{align}
where $\beta_l$ is the $r\times 1$ matrix with entries
$\beta_l^1,\dots,\beta_l^r$.

Consider the $r \times 1$ matrix of derivatives $\partial_{x^l}(H^\top \beta_j)$ 
for fixed  $j,l=1,\dots,k$. We consider below the product of the $2(n-k) \times
r$ matrix
$\left((\widetilde{X}_{1},\tilde\alpha^1),\dots,(\widetilde{X}_r,\tilde\alpha^r)\right)
B$ with the $r \times 1$ matrix 
$\partial_{x^l}(H^\top \beta_j)$. We need the conclusion of Step 2: 
$\left((\widetilde{X}_{1},\widetilde{\alpha}^1),\dots,(\widetilde{X}_{r},\widetilde{\alpha}^r)\right)
B$ does not depend on $x^1, \ldots, x^k$, that is,
\begin{align}\label{zero_bracket_B}
0 &= \partial_{x^l}
\left(\left((\widetilde{X}_{1},\tilde\alpha_1),\dots,(\widetilde{X}_{r},\tilde\alpha_r)\right)B\right)
=\left[(\partial_{x^l},0),\left((\widetilde{X}_{1},\tilde\alpha_1),\dots,(\widetilde{X}_{r},\tilde\alpha_r)
\right)B\right] . \end{align} 
Therefore, the definition \eqref{defb} of the matrix $B$, \eqref{zero_bracket_B}, and the Leibniz rule yield
\begin{align}
\label{Hbeta_intermediary}
&\left((\widetilde{X}_{1},\widetilde{\alpha}_1),\dots,(\widetilde{X}_{r},\widetilde{\alpha}_r)\right)
B \,\partial_{x^l}\left(H^\top\beta_j\right)\\
& \quad =\left[(\partial_{x^l},0),
\left((\widetilde{X}_{1},\tilde\alpha_1),\dots,(\widetilde{X}_{r},\tilde\alpha_r)\right)B\right]
\,\left(H^\top\beta_j\right)
+\left((\widetilde{X}_{1},\tilde\alpha_1),\dots,(\widetilde{X}_{r},\tilde\alpha_r)\right)B\,
\partial_{x^l}\left(H^\top\beta_j\right) \nonumber\\
& \quad =\left[(\partial_{x^l},0),\left((\widetilde{X}_{1},\tilde\alpha_1),\dots,
(\widetilde{X}_{r},\tilde\alpha_r)\right) B H^\top\beta_j\right]\nonumber\\
& \quad =\left[(\partial_{x^l},0),\left((\widetilde{X}_{1},\tilde\alpha_1),\dots,
(\widetilde{X}_{r},\tilde\alpha_r)\right)\beta_j\right] =\left[ (\partial
_{x^l},0), \left[ (\partial _{x^j},0), (\widetilde{X},\tilde\alpha)\right]\right].   \nonumber
\end{align} 
However, 
\begin{align*}
&\left[ (\partial _{x^l},0), \left[ (\partial _{x^j},0), (\widetilde{X},\widetilde\alpha)\right]\right]
=\left[(\partial _{x^l},0), \left(\left[\partial _{x^j},\widetilde X\right],
\boldsymbol{\pounds}_{\partial _{x^j}}\widetilde\alpha\right)\right]\nonumber\\
&\qquad =   \left(\left[\partial _{x^l},\left[\partial _{x^j},\tilde X\right]\right],
\boldsymbol{\pounds}_{\partial _{x^l}}\boldsymbol{\pounds}_{\partial _{x^j}}\tilde\alpha\right)\nonumber\\
& \qquad =\left(- \left[\tilde X,\left[\partial _{x^l},\partial _{x^j}\right]\right]-
  \left[\partial _{x^j},\left[\tilde X,\partial _{x^l}\right]\right] ,
\boldsymbol{\pounds}_{\partial _{x^l}}\boldsymbol{\pounds}_{\partial _{x^j}}\tilde\alpha\right)\nonumber
\end{align*} 
by the Jacobi identity. Since $ [\partial_{x^l},\partial_{x^j}] = 0$  and $
\boldsymbol{\pounds}_{\partial _{x^j}}\boldsymbol{\pounds}_{\partial
  _{x^l}}\widetilde\alpha 
= \boldsymbol{\pounds}_{\partial _{x^l}}\boldsymbol{\pounds}_{\partial
  _{x^j}}\widetilde\alpha$, 
we conclude from  \eqref{Hbeta_intermediary},
\begin{align}
\label{Hbeta}
&\left((\widetilde{X}_{1},\widetilde{\alpha}_1),\dots,(\widetilde{X}_{r},\widetilde{\alpha}_r)\right)
B \,\partial_{x^l}\left(H^\top\beta_j\right)
= \left(\left[\partial _{x^j},\left[\partial _{x^l},\tilde X\right]\right] ,
\boldsymbol{\pounds}_{\partial _{x^j}}\boldsymbol{\pounds}_{\partial _{x^l}}\widetilde\alpha\right) \\
& \qquad = \left[ (\partial _{x^j},0), \left[ (\partial _{x^l},0),
    (\widetilde{X},\tilde\alpha)\right]\right]\nonumber\\
& \qquad =
\left((\widetilde{X}_{1},\tilde\alpha_1),\dots,(\widetilde{X}_{r},\tilde\alpha_r)\right)
B\,\partial_{x^j}\left(H^\top\beta_l\right). \nonumber 
\end{align}

Define the $r\times 1$ matrix with $C^\infty$-entries in the variables $x^1, \ldots, x^n$,  
\begin{align*}
\Pi=(\Pi_{1},\dots,\Pi_r)^\top
:=&-B\left[\int_{0}^{x^k}(H^\top\beta_k)(x^1,\dots,x^{k-1},\tau,x^{k+1},\dots,x^n)d\tau\right.\\
&\quad +\int_{0}^{x^{k-1}}(H^\top\beta_{k-1})(x^1,\dots,x^{k-2},\tau,0,x^{k+1},\dots,x^n)d\tau\\
&\quad \left.+\dots+\int_{0}^{x^{1}}(H^\top\beta_{1})(\tau,0,
\dots,0,x^{k+1},\dots,x^n)d\tau\right] \nonumber \\
=& -(BR)(x^1, \ldots, x^n),
\end{align*}
where $R(x^1, \ldots, x^n)$ is the $r\times 1$ matrix  in the parenthesis.
Then  for $l=1,\ldots,k$, we get
\begin{align}\label{diffZ}
&\left[(\partial_{x^l},0),\left((\widetilde{X},\tilde\alpha)
+\left((\widetilde{X}_{1},\tilde\alpha_1),\dots,(\widetilde{X}_{r},\tilde\alpha_r)\right)
\Pi\right)\right]   \\
& \qquad =  \left[(\partial_{x^l},0),(\widetilde{X},\tilde\alpha) \right] 
+ \left[(\partial_{x^l},0), \left((\widetilde{X}_{1},\tilde\alpha_1),\dots,
(\widetilde{X}_{r},\tilde\alpha_r)\right)\Pi  \right]  \nonumber  \\
& \qquad =  \left[(\partial_{x^l},0),(\widetilde{X},\tilde\alpha) \right] 
- \left[(\partial_{x^l},0), \left((\widetilde{X}_{1},\tilde\alpha_1),\dots,
(\widetilde{X}_{r},\tilde\alpha_r)\right)BR  \right]  \nonumber  \\
& \qquad  = \left[(\partial_{x^l},0),(\widetilde{X},\tilde\alpha) \right] -
\left((\widetilde{X}_{1},\tilde\alpha_1),\dots,(\widetilde{X}_{r},\tilde\alpha_r)\right)
B\left(\partial_{x^l} R\right)  \nonumber 
\end{align} 
by \eqref{I_bracket_function} and using the fact that 
$\left((\widetilde{X}_{1},\tilde\alpha_1),\dots,(\widetilde{X}_{r},\tilde\alpha_r)\right)
B$ is independent of $x^1, \ldots, x^k$  (see the conclusion of Step 2). 
For any $l=1, \ldots, k$, we prove the identity 
\begin{equation}\label{R_id}
\left((\widetilde{X}_{1},\tilde\alpha_1),\dots,(\widetilde{X}_{r},\tilde\alpha_r)\right)
B\partial_{x^l} R
= \left((\widetilde{X}_{1},\tilde\alpha_1),\dots,(\widetilde{X}_{r},\tilde\alpha_r)\right) \beta_l
\end{equation} 
Indeed, since we can move freely 
$
\left((\widetilde{X}_{1},\tilde\alpha_1),\dots,(\widetilde{X}_{r},\tilde\alpha_r)\right)
B $ (a matrix depending smoothly on $x^{k+1}, \ldots, x^n$) 
in and out of the integral signs in the computation below, we get
\begin{align*}
&\left((\widetilde{X}_{1},\tilde\alpha_1),\dots,
(\widetilde{X}_{r},\tilde\alpha_r)\right) B\partial_{x^l} R \\
&= \left((\widetilde{X}_{1},\tilde\alpha_1),\dots,(\widetilde{X}_{r},\tilde\alpha_r)
\right) B\partial_{x^l} \left[\int_{0}^{x^k}(H^\top\beta_k)(x^1,\dots,x^{k-1},\tau,
x^{k+1},\dots,x^n)d\tau\right.\\
& \quad +\int_{0}^{x^{k-1}}(H^\top\beta_{k-1})(x^1,\dots,x^{k-2},\tau,0,x^{k+1},
\dots,x^n)d\tau\\
&\quad \left.+\dots+\int_{0}^{x^{1}}(H^\top\beta_{1})(\tau,0,
\dots,0,x^{k+1},\ldots,x^n)d\tau \right]\nonumber \\
& = \int_{0}^{x^k} \left((\widetilde{X}_{1},\tilde\alpha_1),\dots,
(\widetilde{X}_{r},\tilde\alpha_r)\right) B 
\partial_{x^l} (H^\top\beta_k)(x^1,\dots,x^{k-1},\tau,
x^{k+1},\dots,x^n)d\tau \\
& \quad +\int_{0}^{x^{k-1}}
\left((\widetilde{X}_{1},\tilde\alpha_1),\dots,
(\widetilde{X}_{r},\tilde\alpha_r)\right) B 
\partial_{x^l} (H^\top\beta_{k-1})(x^1,\dots,x^{k-2},\tau,0,x^{k+1},\dots,x^n)d\tau\\
&\quad +\dots+\int_{0}^{x^{l+1}} \left((\widetilde{X}_{1},\tilde\alpha_1),\dots,(\widetilde{X}_{r},\tilde\alpha_r)\right) B \partial_{x^l} (H^\top\beta_{l+1})(x^1,\dots,x^{l},\tau,0,\ldots,0,x^{k+1},\dots,x^n)d\tau\\
&\quad +\partial_{x^l}\int_{0}^{x^{l}} \left((\widetilde{X}_{1},\tilde\alpha_1),\dots,(\widetilde{X}_{r},\tilde\alpha_r)\right) B(H^\top\beta_{l})(x^1,\dots,x^{l-1},\tau,0,\ldots,0,x^{k+1},\dots,x^n)d\tau\\
&\quad +\int_{0}^{x^{l-1}} \left((\widetilde{X}_{1},\tilde\alpha_1),\dots,(\widetilde{X}_{r},\tilde\alpha_r)\right) B \partial_{x^l} (H^\top\beta_{l-1})(x^1,\dots,x^{l-2},\tau,0,\ldots,0,x^{k+1},\dots,x^n)d\tau\\
&\quad +\dots +\int_{0}^{x^{1}} \left((\widetilde{X}_{1},\tilde\alpha_1),\dots,(\widetilde{X}_{r},\tilde\alpha_r)\right) B \partial_{x^l} (H^\top\beta_{1})(\tau,0,
\dots,0,x^{k+1},\ldots,x^n)d\tau. \nonumber 
\end{align*} 
Since $(H^\top\beta_{m+1})(x^1,\ldots,x^m,\tau,0,\ldots,0,x^{k+1},\ldots,x^n)$ 
doesn't depend on $x^l$ for $m< l-1$, the last $l-1$ integrals in this
expression vanish. Using \eqref{Hbeta} in the first $k-l$ integrals we get

\begin{align*}
&\int_{0}^{x^k} \left((\widetilde{X}_{1},\tilde\alpha_1),\dots,(\widetilde{X}_{r},\tilde\alpha_r)\right) B \partial_{x^k} (H^\top\beta_l)(x^1,\dots,x^{k-1},\tau, x^{k+1},\dots,x^n)d\tau \\
& \quad +\int_{0}^{x^{k-1}}
\left((\widetilde{X}_{1},\tilde\alpha_1),\dots,(\widetilde{X}_{r},\tilde\alpha_r)\right) B \partial_{x^{k-1}}(H^\top\beta_{l})(x^1,\dots,x^{k-2},\tau,0,x^{k+1},\dots,x^n)d\tau\\
&\quad +\dots +\int_{0}^{x^{l+1}} \left((\widetilde{X}_{1},\tilde\alpha_1),\dots,(\widetilde{X}_{r},\tilde\alpha_r)\right) B \partial_{x^{l+1}} (H^\top\beta_{l})(x^1,\dots,x^{l},\tau,0,\ldots,0,x^{k+1},\dots,x^n)d\tau\\
&\quad +\partial_{x^l}\int_{0}^{x^{l}} \left((\widetilde{X}_{1},\tilde\alpha_1),\dots,(\widetilde{X}_{r},\tilde\alpha_r)\right) B(H^\top\beta_{l})(x^1,\dots,x^{l-1},\tau,0,\ldots,0,x^{k+1},\dots,x^n)d\tau\\
& = \left((\widetilde{X}_{1},\tilde\alpha_1),\dots,(\widetilde{X}_{r},\tilde\alpha_r)\right) B\\
&\quad \times\Biggr[\Bigl((H^\top \beta _l )(x^1, \ldots, x^k, \ldots, x^n) 
- (H^\top \beta _l )(x^1, \ldots, x^{k-1}, 0, x^{k+1}, \ldots, x^n) \Bigr)\\
&\quad +\Bigl((H^\top \beta _l )(x^1, \ldots, x^{k-1}, 0, x^{k+1}, \ldots, x^n)
- (H^\top \beta _l )(x^1, \ldots, x^{k-2}, 0, 0, x^{k+1}, \ldots, x^n)\Bigr)\\
&\quad +\ldots+\Bigl((H^\top \beta _l )(x^1, \ldots, x^{l+1},0,\ldots, 0, x^{k+1}, \ldots, x^n)
- (H^\top \beta _l )(x^1, \ldots, x^{l}, 0,\ldots, 0, x^{k+1}, \ldots, x^n)\Bigr)\\
&\quad   + (H^\top \beta _l )(x^1,\ldots,x^l,0,\ldots, 0, x^{k+1}, \ldots, x^n)\Biggr]\\
&\quad = \left((\widetilde{X}_{1},\tilde\alpha_1),\dots,(\widetilde{X}_{r},\tilde\alpha_r)\right) B (x^1, \ldots,  x^n) H^\top(x^1, \ldots,  x^n) \beta _l(x^1, \ldots,  x^n)\\
&\quad = \left((\widetilde{X}_{1},\tilde\alpha_1),\dots,(\widetilde{X}_{r},\tilde\alpha_r)\right) \beta _l(x^1, \ldots,  x^n)
\end{align*} 
by \eqref{defb}. This proves \eqref{R_id}. 
 
From \eqref{R_id}, \eqref{diffZ},  and \eqref{bracket_tilde_l} we conclude
\begin{align}\label{diffZ_next}
&\left[(\partial_{x^l},0),\left((\widetilde{X},\widetilde\alpha)+\left((\widetilde{X}_{1},\widetilde\alpha_1),\dots,(\widetilde{X}_{r},\widetilde\alpha_r)\right)
\Pi\right)\right]   \\
& \qquad = \left((\widetilde{X}_{1},\widetilde\alpha_1),\dots,(\widetilde{X}_{r},\widetilde\alpha_r)\right)\beta_l-\left((\widetilde{X}_{1},\widetilde\alpha_1),\dots,(\widetilde{X}_{r},\widetilde\alpha_r)\right)\beta_l=0. \nonumber 
 \end{align}
 
This identity suggests that the required section $d=(Z, \gamma ) \in  \Gamma(
\mathcal{D})$ satisfying the third condition in the statement of the
proposition is
\begin{align*}
d=(Z,\gamma):=\left((X_{1},\alpha_1),\ldots,(X_{r},\alpha_r)\right)\Pi
=\sum_{l=1}^r(X_l,\alpha_l)\Pi_l\in \Gamma(\mathcal{D}).
\end{align*}
We have
\begin{align*}
(Z, \gamma ) &=\left((\bar X_{1},0),\ldots,(\bar X_{r},0)\right)\Pi+\left((\widetilde X_{1},\widetilde\alpha_1),\ldots,(\widetilde X_{r},\widetilde\alpha_r)\right)\Pi\\
&=\sum_{l=1}^r(\bar X_l,0)\Pi_l+\sum_{l=1}^r(\widetilde X_l,\widetilde\alpha_l)\Pi_l,
\end{align*}
where $\bar{X}:=\sum_{j=1}^ka^j\partial_{x^j}$ for $i=1,\dots,r$. Set in the
same manner $\bar{X}_i=\sum_{j=1}^kX_i^j\partial_{x^j}\in\Gamma(\I)$ for
$i=1,\dots,r$, and  verify (iii) in the statement of the theorem. For any $l = 1, \ldots, k$, since 
$ \partial_{x^1}, \ldots, \partial_{x^k} $
is a basis of the space of sections of $\mathcal{I}$ over $U$, we get
 \begin{align*} 
\left[ (X+Z,\alpha+\gamma), (\partial _{x^l},0)\right]=&\left[ (\bar X,0)+\sum_{k=1}^r(\bar
  X_k,0)\Pi_k, (\partial _{x^l},0)\right]+\left[ (\tilde X,\tilde\alpha)+\sum_{k=1}^r(\tilde
  X_k,\tilde\alpha_k)\Pi_k, (\partial _{x^l},0)\right]\\
\overset{\eqref{diffZ_next}}=&\left[ (\bar X,0)+\sum_{k=1}^r(\bar
  X_k,0)\Pi_k, (\partial _{x^l},0)\right]+0\in  \Gamma (\Theta)
\end{align*}
since, by construction, $\bar{X}, \bar{X}_i \in  \Gamma( \mathcal{I})$,
$i=1,\dots,r$, so that $\bar{X}+\sum_{k=1}^r\bar{X}_k\Pi_k\in  \Gamma ( \I )$.
Since $( \partial_{x^1}, 0), \ldots,$ $( \partial_{x^k}, 0) $ span 
the distribution $ \Theta$ over $U$, we conclude that 
$[(X+Z,\alpha+\gamma),\Gamma(\Theta)]\subseteq\Gamma(\Theta)$ on $U$ 
and  (iii) in the statement is proved.

\subsection{Example}
Let $(G\rr M, \tg,\s,\m, \epsilon)$ be a $\tg$-connected Lie groupoid with Lie algebroid $A\to M$.
The Pontryagin bundle $TG\oplus T^*G$
inherits the structure of a Lie groupoid
over the vector bundle $TM\oplus A^*$ (see for instance \cite{Jotz10b}). The source and target maps
of this groupoid are written $\mathbb T\s$ and $\mathbb T\tg$.
Consider a smooth subdistribution $\mathcal D$ of  the kernel $\ker\mathbb 
T\s=T^\s G\oplus (T^\tg G)^\circ$, where $T^\s G$ is the 
involutive subbundle tangent to the $\s$-fibers and 
$T^\tg G\subseteq TG$ is the involutive subbundle
tangent to the $\tg$-fibers. 

For any left-invariant section 
$X^l$ of $T^\tg G$ and for any 
right-invariant section $(Y^r,\tg^*\alpha)$ of 
$\ker\mathbb{T}\s$, we find that 
$\left[\left(X^l,0\right), \left(Y^r,\tg^*\alpha\right)\right]
=0$ (see for instance \cite{Jotz10b}). Since these sections 
span $T^\tg G$ and $\ker\mathbb{T}\s$, respectively,
we find that
$[\Gamma(\ker\mathbb T\s), \Gamma(T^\tg G\oplus\{0\})]
\subseteq \Gamma(\ker\mathbb T\s+T^\tg G\oplus\{0\})$.

Choose $g\in G$ and assume that there exist smooth 
sections $d_1,\ldots,d_r$ of $\mathcal D$, defined on a neighborhood
$U$ of $g$ in $G$ and such that
\[\mathcal D(g')=\erz_\R\{d_1(g'), \ldots, d_r(g')\}\]
for all $g'\in U$ and 
\begin{equation}\label{ass1}
[d_i,\Gamma\left(T^\tg G\oplus\{0\}\right)]\subseteq\erz_{C^\infty(U)}\{d_1,\ldots,d_r\}
+\Gamma\left(T^\tg G\oplus\{0\}\right)
\end{equation}
for $i=1,\ldots,r$.

Then, using Theorem
\ref{prop-ch-ta}, we conclude that, on $U$,  $\mathcal D$
is spanned pointwise by sections of the type $(X,\tg^*\alpha)\in\Gamma(\ker\mathbb T\s)$
with $X\sim_\tg\bar X$ for some $\bar X\in\mx(M)$ and $\alpha\in\Omega^1(M)$, 
i.e., such that $\mathbb T\tg\circ (X,\tg^*\alpha)$ is constant on $\tg$-fibers.

Furthermore, if $p\in\Gamma(\mathsf P)$ is defined on $U$ and 
such that
$[p, \Gamma(T^\tg G\oplus\{0\})]
\subseteq \erz_{C^\infty(U)}\{d_1,\ldots,d_r\}+\Gamma\left(T^\tg G\oplus\{0\}\right)$,
then there exists $k\in\Gamma(\mathcal D)$
such that $p+k$ is $T^\tg G$-invariant
and $\mathbb T\s\circ(p+k)=\mathbb T\s\circ p$.

\medskip

For another example from control theory, see \cite{ChTa89}.

\bibliographystyle{amsalpha}

\def\cprime{$'$} \def\polhk#1{\setbox0=\hbox{#1}{\ooalign{\hidewidth
  \lower1.5ex\hbox{`}\hidewidth\crcr\unhbox0}}} \def\cprime{$'$}
  \def\cprime{$'$} \def\cprime{$'$} \def\cprime{$'$} \def\cprime{$'$}
  \def\cprime{$'$} \def\cprime{$'$}
  \def\polhk#1{\setbox0=\hbox{#1}{\ooalign{\hidewidth
  \lower1.5ex\hbox{`}\hidewidth\crcr\unhbox0}}}
  \def\polhk#1{\setbox0=\hbox{#1}{\ooalign{\hidewidth
  \lower1.5ex\hbox{`}\hidewidth\crcr\unhbox0}}}
  \def\polhk#1{\setbox0=\hbox{#1}{\ooalign{\hidewidth
  \lower1.5ex\hbox{`}\hidewidth\crcr\unhbox0}}}
  \def\polhk#1{\setbox0=\hbox{#1}{\ooalign{\hidewidth
  \lower1.5ex\hbox{`}\hidewidth\crcr\unhbox0}}} \def\cprime{$'$}
  \def\polhk#1{\setbox0=\hbox{#1}{\ooalign{\hidewidth
  \lower1.5ex\hbox{`}\hidewidth\crcr\unhbox0}}}
  \def\polhk#1{\setbox0=\hbox{#1}{\ooalign{\hidewidth
  \lower1.5ex\hbox{`}\hidewidth\crcr\unhbox0}}}
  \def\polhk#1{\setbox0=\hbox{#1}{\ooalign{\hidewidth
  \lower1.5ex\hbox{`}\hidewidth\crcr\unhbox0}}}
  \def\polhk#1{\setbox0=\hbox{#1}{\ooalign{\hidewidth
  \lower1.5ex\hbox{`}\hidewidth\crcr\unhbox0}}}
\providecommand{\bysame}{\leavevmode\hbox to3em{\hrulefill}\thinspace}
\providecommand{\MR}{\relax\ifhmode\unskip\space\fi MR }
\providecommand{\MRhref}[2]{%
  \href{http://www.ams.org/mathscinet-getitem?mr=#1}{#2}
}
\providecommand{\href}[2]{#2}

\end{document}